\documentclass{article}

\usepackage{amssymb,amsmath,fullpage,amsthm,amsfonts, verbatim, tikz, bbm}
\usepackage{url}
\usepackage{algorithm}
\usepackage{algpseudocode}
\usepackage{xcolor}
\definecolor{linkblue}{HTML}{003d73}
\definecolor{linkgreen}{HTML}{006161}
\definecolor{linkred}{HTML}{a11950}
\usepackage{hyperref}
\hypersetup{
	pdftitle={Expected Distances on Manifolds of Partially Oriented Flags},
	pdfauthor={Brenden Balch, Chris Peterson, and Clayton Shonkwiler},
	pdfsubject={differential geometry},
	pdfkeywords={flag manifolds, random geometry, data models},
	colorlinks=true,
	linkcolor=linkblue,
	citecolor=linkgreen,
	urlcolor=linkred
}
\usepackage{authblk}
\usepackage{ytableau}
\usepackage[nameinlink]{cleveref}
\usepackage{bigints}

\usetikzlibrary{arrows}
\usetikzlibrary{positioning}
\usetikzlibrary{cd}

\newtheorem{theorem}{Theorem}
\newtheorem{prop}[theorem]{Proposition}

\theoremstyle{definition}
\newtheorem{ex}{Example}[section]
\newtheorem{definition}{Definition}
\newcommand{\R}{\mathbb{R}}
\newcommand{\Q}{\mathbb{Q}}
\newcommand{\C}{\mathbb{C}}
\newcommand{\N}{\mathbb{N}}

\newcommand{\Sph}{\mathbb{S}}
\newcommand{\vol}{\operatorname{Vol}}

\newcommand{\Fl}{F\ell}
\newcommand{\Gr}{\operatorname{Gr}}
\newcommand{\St}{\operatorname{St}}
\newcommand{\dVol}{\operatorname{dVol}}
\newcommand{\dArea}{\operatorname{dArea}}
\newcommand{\Area}{\operatorname{Area}}
\newcommand{\E}{\mathbb{E}}

\title{Expected Distances on Manifolds of Partially Oriented Flags}
\author{Brenden Balch}
\author{Chris Peterson}
\author{Clayton Shonkwiler}
\affil{Department of Mathematics, Colorado State University, Fort Collins, CO}
\date{}

\begin{document}
	
\maketitle

\begin{abstract}
	Flag manifolds are generalizations of projective spaces and other Grassmannians: they parametrize flags, which are nested sequences of subspaces in a given vector space. These are important objects in algebraic and differential geometry, but are also increasingly being used in data science, where many types of data are properly understood as subspaces rather than vectors. In this paper we discuss partially oriented flag manifolds, which parametrize flags in which some of the subspaces may be endowed with an orientation. We compute the expected distance between random points on some low-dimensional examples, which we view as a statistical baseline against which to compare the distances between particular partially oriented flags coming from geometry or data.
	\end{abstract}

%----------------------------------------------------------------------------------------------------------------------------------
\section{Introduction}
Let \(V\) be an $n$-dimensional vector space. A \textbf{flag} in \(V\) is a nested sequence of subspaces ${V_1 \subset V_2 \subset \dots \subset V_k =V}$. If $d_i$ is the dimension of $V_i$, then associated to the flag is an increasing sequence ${d_1<d_2<\dots < d_k=n}$. We call the sequence $(d_1,d_2, \dots, d_k)$ the \textbf{signature} of the flag. For this paper, we will make the assumption that $V$ is a real vector space and the flag manifold $\Fl(d_1,d_2, \dots, d_k)$ will be the real manifold whose points parametrize all flags of signature $(d_1,d_2, \dots, d_k)$. Flag manifolds are natural generalizations of Grassmann manifolds and have found applications in image analysis~\cite{Nishimori:2006ke} and face, pose, and action recognition~\cite{Draper:2014gu, Marrinan:2015cj}; as Ye, Wong, and Lim point out~\cite{Ye:2019vy}, they are also implicit in many tasks in numerical and statistical analysis, ranging from mesh refinement to multiresolution analysis to canonical correlation analysis. The key feature of flag manifolds is they locally look the same at each point, i.e. flag manifolds are homogeneous spaces. Homogeneous spaces naturally admit nice coordinates which allow for tangible calculations.

In this paper we generalize the notion of flag manifold, defining \textbf{partially oriented flag manifolds}. These manifolds are slight generalizations of the partially oriented flag manifolds introduced by Lam~\cite{Lam:1975kj} and studied by Sankaran and Zvengrowski~\cite{Sankaran:1987wh,Sankaran:1997cj}, among others. Flag manifolds are natural models for data which is properly represented as (nested) subspaces, rather than as individual vectors, and partially oriented flag manifolds should prove useful as models for subspace data when some or all of the nested subspaces are equipped with an orientation. These spaces are rather simple to define as homogeneous spaces, but do not seem to have appeared previously in the data science literature.

Given two data, represented as points in some Riemannian manifold or more general metric space, it is natural to use the distance between the points as a measure of similarity. However, the raw distance is meaningless without further context to know whether the distance is, say, much smaller than expected. A reasonable statistical baseline for such a comparison is the expected distance between two random points in the space. With an eye towards applications of partially oriented flag manifolds to data problems, we determine the expected geodesic distance between random points in all manifolds of partially oriented flags in \(\R^3\); compare to Absil, Edelman, and Koev's work using a slightly different notion of distance in Grassmannians~\cite{edelman2}.

More precisely, recall that \(\Fl(1,2,3)\) is the manifold of all \((1,2,3)\)-flags in \(\R^3\); that is, each point represents a line inside a plane inside \(\R^3\). We will shortly introduce some more involved notation for this space, but for the moment let \(M\) be the manifold of such flags in which the line has a preferred orientation, but the plane and \(\R^3\) do not. Then the most interesting expectations we compute are:

\begin{theorem}\label{thm:main theorem}
	The expected distance between two random points in \(M\) is
	\[
		\E[d;M] = 1 + \frac{\pi}{4}.
	\]
	The expected distance between random \((1,2,3)\)-flags is
	\begin{multline*}
		\E[d;\Fl(1,2,3)] = \frac{3\pi}{2} + \frac{96}{\pi^2} \bigintsss_0^{\pi/4} \! \left[\vphantom{\frac{\arctan^2 \left(\sqrt{1+\sec^2 \varphi_3}\right)}{\sqrt{1+\sec^2\varphi_3}}} \arctan\left( \tan^2 \left(\frac{\arctan(\sec \varphi_3)}{2}\right)\right)\right. \\
		\left. - \frac{\arctan^2 \left(\sqrt{1+\sec^2 \varphi_3}\right)}{\sqrt{1+\sec^2\varphi_3}}\right]\! d\varphi_3 \approx 1.31172503.
	\end{multline*}
\end{theorem}

We define the partially oriented flag manifolds in \Cref{sec:flags} and discuss a number of examples and relationships to more familiar manifolds, including Grassmannians, Stiefel manifolds, and classical flag manifolds. In \Cref{sec:so3} we compute the expected distance between two random points in \(SO(3)\) both analytically and numerically using Monte Carlo integration. While this computation is not new, it serves as a starting point and template for determining the expected distances between random points in all the manifolds of (partially oriented) flags on \(\R^3\), which is the focus of \Cref{sec:expected distances}.

%----------------------------------------------------------------------------------------------------------------------------------
\section{Flags and Partially Oriented Flags}\label{sec:flags}
Let $\mathbb N$ denote the set of all positive integers and let $n\in \mathbb N$. An {\bf ordered partition of the integer  \(n\)} is the tuple \(\lambda = (\lambda_1, \dots, \lambda_k)\) with \(\lambda_j \in \N\) such that \(\sum_{j=1}^k\lambda_j = n\).\footnote{Contrast this with the usual definition of a partition, in which the order does not matter, and hence one can assume \(\lambda_1 \geq \lambda_2, \dots, \geq \lambda_k\).} In this case, we say the {\bf size} of \(\lambda\) is \(k\). Notice that for an ordered partition of size \(k\), the symmetric group \(S_k\) acts by permuting elements; {\it i.e.} if \(\sigma \in S_k\) we define \(\sigma \cdot \lambda = (\lambda_{\sigma(1)}, \dots, \lambda_{\sigma(k)})\). Now suppose that \(I_\lambda = \{1,\dots , k\}\) is the collection of indices of the partition \(\lambda\). Then a {\bf partition of the set \(I_\lambda\)} is a collection of disjoint nonempty subsets of \(I_\lambda\) whose union is all of \(I_\lambda\). We call the partition \(\{I_\lambda\}\) the {\bf trivial partition}, and the partition \( \{ \{1\}, \{2\} , \dots , \{k\}\}\) the {\bf complete partition}. All other partitions of \(I_\lambda\) will be called {\bf proper partitions}.

Let $O(r)$ denote the group of real $r\times r$ orthogonal matrices. For any fixed ordered partition \(\lambda\) of \(n\), define the group \(G_\lambda = \prod_{j \in I_\lambda}O(\lambda_j)\). An easy exercise in group theory shows that if \(\sigma \in S_k\), then \(G_\lambda \cong G_{\sigma \cdot \lambda}\). Since elements of \(G_\lambda\) are of the form \(\bigoplus_{j \in I_\lambda} A_j\) where \(A_j \in O(\lambda_j)\), we can interpret \(G_\lambda\) as a block-diagonal subgroup of \(O(n)\). 

Now suppose that \(I_\alpha = \{{i_1}, \dots, {i_m}\}\) is a subset of \(I_\lambda\). We define a {\bf special orthogonal block of \(G_\alpha\)} to be the set 
\[
S(O(\lambda_{i_1}) \times \cdots \times O(\lambda_{i_m})) = \{A \in G_\alpha | \det(A) =1\}.
\]
For brevity of notation, we will denote this set as \(SG_\alpha\). If \(P\) is a partition of \(I_\lambda\), define a {\bf special orthogonal partition of \(G_\lambda\)} to be the group \(SG_\lambda^P = \prod_{I_\alpha \in P}SG_{\alpha}\). If the partition of \(I_\lambda\) is trivial, complete, or proper, we call the corresponding special orthogonal partition the same. \(SG_\lambda^P\) is a block-diagonal subgroup of \(SO(n)\) which is finite if and only if \(\lambda = (1,1,\dots1)\). Since the \(SG_\lambda^P\) are subgroups of \(SO(n)\), the quotients will be homogeneous spaces which are the central object of this paper.

\begin{definition}
Suppose that \(\lambda\) is any ordered partition of the integer \(n\). Then if \(P\) is a partition of \(I_\lambda\), we may define the orbit space \(\Fl(\lambda;P) = SO(n)/SG_\lambda^P\). If the partition \(P\) is trivial, then \(\Fl(\lambda;P)\) is called a {\bf flag manifold}. If \(P\) is a complete partition, we say that \(\Fl(\lambda;P)\) is an {\bf oriented flag manifold}. Finally, if the partition \(P\) is proper, then \(\Fl(\lambda;P)\) is called a {\bf partially oriented flag manifold}.
\end{definition}

\begin{ex}
\label{grass}
If \(\lambda = (k,n-k)\) and \(P_T = \{\{1, 2\}\} \) is the trivial partition, then
\[
	\Fl(\lambda;P_T) = SO(n)/S(O(k) \times O(n-k))\cong \Gr(k,n),
\] 
the Grassmannian of $k$-dimensional linear subspaces of $\R^n$. If instead we take the complete partition \(P_C = \{\{1\},\{2\}\}\), then 
\[
	\Fl(\lambda;P_C) = SO(n)/(SO(k) \times SO(n-k)) \cong \widetilde{\Gr}(k,n),
\]
the oriented Grassmannian of oriented \(k\)-dimensional subspaces of \(\R^n\), which double covers \(\Gr(k,n)\).
\end{ex}

\begin{ex}
If $1\leq d_1 < d_2 < \dots < d_k = n$, recall that the flag manifold \(\Fl(d_1, \dots d_k)\) consists of nested subspaces
\[
V_1 \subset V_2 \subset \cdots \subset V_k = \R^n
\]
where \(\dim(V_m) = d_m\).

If \(\lambda = (\lambda_1, \dots \lambda_k)\) is an ordered partition of the integer \(n\), define \(d_m = \sum_{j=1}^m\lambda_j\) for each \(1\leq m \leq k\). If \(P_T\) is the trivial partition of \(I_\lambda\), then 
\[
\Fl(\lambda;P_T) = SO(n)/S(O(\lambda_1)  \times \cdots \times O(\lambda_k)) \cong \Fl(d_1, \dots d_k)
\]
really is a flag manifold. Notice that if \(\sigma \in S_k\), then \(\Fl(\sigma \cdot \lambda;P_T)\) yields a different space parametrizing entirely different geometric objects which is nonetheless diffeomorphic to \(\Fl(\lambda;P_T)\). This helps explain our use of ordered partitions of \(n\) rather than just partitions: we are thinking of the \(\Fl(\lambda;P)\) spaces as data models, with points corresponding to actual data, so we want to avoid the complication of applying a diffeomorphism before we can think of our data as points in a (partially oriented) flag manifold. 
\end{ex}

\begin{ex}
Take \(\lambda = (n-k,1,1,\dots,1)\) and suppose that \(P_C\) is the complete partition of \(I_\lambda\). Then 
\[
\Fl(\lambda;P_C) = SO(n)/(SO(n-k) \times SO(1) \times \cdots \times SO(1)) \cong SO(n)/SO(n-k) \cong \St(k,n),
\] 
the Stiefel manifold of orthonormal \(k\)-frames in \(\R^n\). 
\end{ex}

\begin{ex}
\label{fin}
In this example, we'll consider all special orthogonal partitions when \(n=3\) and \(\lambda = (1,1,1)\). If $P_C=\{\{1\},\{2\},\{3\}\}$ is the complete partition of \(I_\lambda\), then \(SG_\lambda^{P_C} = SO(1) \times SO(1) \times SO(1)\) is the trivial group, so \(\Fl(\lambda;P_C) \cong SO(3)\). At the other extreme, if $P_T=\{\{1,2,3\}\}$ is the trivial partition, then \(\Fl(\lambda;P_T) \cong Fl(1,2,3))\) and  \(SG_\lambda^{P_T} = S(O(1) \times O(1) \times O(1)) \) is the copy of the Klein 4-group
\[SG_\lambda^{P_T} =  \left\{
\begin{pmatrix}
1 & 0 & 0\\
0 & 1 & 0\\
0 & 0 & 1\\
\end{pmatrix},
\begin{pmatrix}
-1 & 0 & 0\\
0 & -1 & 0\\
0 & 0 & 1\\
\end{pmatrix},
\begin{pmatrix}
1 & 0 & 0\\
0 & -1 & 0\\
0 & 0 & -1\\
\end{pmatrix},
\begin{pmatrix}
-1 & 0 & 0\\
0 & 1 & 0\\
0 & 0 & -1\\
\end{pmatrix}
\right\}.\]

If \(P_1 = \{\{1\},\{2,3\}\}\), \(P_2 = \{\{2\},\{1,3\}\}\), and \(P_3 =\{\{3\},\{1,2\}\} \) are proper partitions, then the corresponding groups are
\[
SG_\lambda^{P_1} = 
\left\{
\begin{pmatrix}
1 & 0 & 0\\
0 & 1 & 0\\
0 & 0 & 1\\
\end{pmatrix},
\begin{pmatrix}
1 & 0 & 0\\
0 & -1 & 0\\
0 & 0 & -1\\
\end{pmatrix}
\right\}
\]

\[
SG_\lambda^{P_2} = 
\left\{
\begin{pmatrix}
1 & 0 & 0\\
0 & 1 & 0\\
0 & 0 & 1\\
\end{pmatrix},
\begin{pmatrix}
-1 & 0 & 0\\
0 & 1 & 0\\
0 & 0 & -1\\
\end{pmatrix}
\right\}
\]

\[
SG_\lambda^{P_3} = 
\left\{
\begin{pmatrix}
1 & 0 & 0\\
0 & 1 & 0\\
0 & 0 & 1\\
\end{pmatrix},
\begin{pmatrix}
-1 & 0 & 0\\
0 & -1 & 0\\
0 & 0 & 1\\
\end{pmatrix}
\right\}.
\]
Clearly each \(SG_\lambda^{P_i}\) is a subgroup of \(SG_\lambda^{P_T}\) of index 2, so in the quotients we have natural double covers captured by the tower in Figure~\ref{fig:tower}.
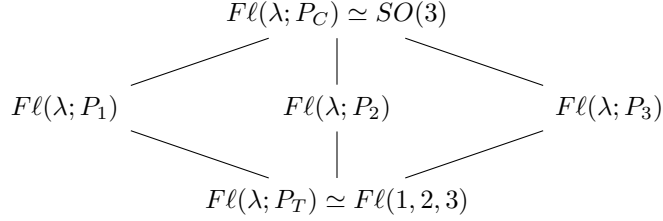
\begin{figure}
\centering
\begin{tikzcd}
	& \Fl(\lambda; P_C) \simeq SO(3) \ar[dl, dash] \ar[d, dash] \ar[dr, dash] & \\
	\Fl(\lambda; P_1) \ar[dr, dash] & \Fl(\lambda; P_2) \ar[d, dash] & \Fl(\lambda; P_3) \ar[dl, dash] \\
	& \Fl(\lambda; P_T) \simeq \Fl(1,2,3) & 
\end{tikzcd}
	
\caption{Tower of flags for \(\lambda = (1,1,1)\).}	
\label{fig:tower}
\end{figure}
\end{ex}

The double coverings in the previous example are special cases of a more general phenomenon. Given an ordered partition \(\lambda\) of \(n\) and any partition \(P\) of \(I_\lambda\), suppose that \(P'\) is a refinement of \(P\). Then the number of sets contained in \(P'\) and \(P\) differ by an integer \(m\). It follows that 
\(|SG_\lambda^P:SG_\lambda^{P'}| = 2^m\), and hence we have the proposition:
\begin{prop}
Suppose that \(\lambda\) is an ordered partition of \(n\) and \(P\) is a set partition of \(I_\lambda\). Given any refinement \(P'\) of \(P\) with \(m = |P'| - |P|\), \(\Fl(\lambda;P')\) is a \(2^m\)-cover of \(\Fl(\lambda;P)\).
\end{prop}

\begin{ex}
In light of the proposition and Example \ref{grass}, we re-verify that \(\widetilde{\Gr}(k,n)\) is a double cover of \(\Gr(k,n)\).
\end{ex}

\begin{figure}[htbp]
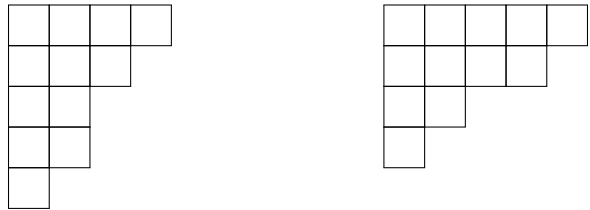

	\centering
		\ydiagram{4,3,2,2,1} \qquad \qquad \qquad \qquad \ydiagram{5,4,2,1}
	\caption{On the left is the Young diagram corresponding to the partition \(\lambda = (4,3,2,2,1)\). The conjugate Young diagram is on the right, with corresponding partition \(\bar{\lambda} = (5,4,2,1)\).}
	\label{fig:Young diagrams}
\end{figure}

Since \(\Fl(\lambda;P)\) is defined by quotients of orthogonal groups, we can find the total volume of  \(\Fl(\lambda;P)\). If \(\lambda\) is of length \(k\) and \(\sigma \in S_k\), we have \({\vol(\Fl(\lambda;P)) = \vol(\Fl(\sigma \cdot \lambda;\sigma \cdot P))}\), so, for the purposes of computing volume, we may assume that \(\lambda\) is in decreasing order.
To simplify notation, let $\Sph^r$ denote the unit $r$-sphere and set \(V_{i} = \vol(\Sph^{i-1})\) for \({i \in \{1,2,\dots n\}}\); then
\[
	\vol(O(m)) = \prod_{i=1}^m V_i,
\]
so we see:

\begin{prop}\label{prop:flag volumes}
	With notation as above,
	\[
		\vol(\Fl(\lambda;P)) =2^{|P|-1}  \frac{V_1 V_2 \cdots V_n}{\prod_{i=1}^{\lambda_1} V_i \cdots \prod_{i=1}^{\lambda_k} V_i} = 2^{|P|-1} \frac{V^{\mathbbm{1}}}{V^{\bar{\lambda}}}=2^{|P|-1}V^{\mathbbm{1}-\bar{\lambda}}.
	\]
	In this expression, \(\bar{\lambda}\) denotes the conjugate of the partition \(\lambda\) given by taking the Young diagram corresponding to \(\lambda\), reflecting across the main diagonal, and then taking the partition corresponding to this new diagram; see \Cref{fig:Young diagrams} for an example.\footnote{It is in converting the partition \(\lambda\) into a Young diagram that it is important we assume \(\lambda\) is in decreasing order.} Also, we use the notation \(\mathbbm{1} := (1, \dots , 1)\) and \(x^\mu := x_1^{\mu_1}x_2^{\mu_2} \cdots x_n^{\mu_n}\).
\end{prop}

%---------------------------------------------------------------------------------------------------------------------------------
%---------------------------------------------------------------------------------------------------------------------------------
\section{\(\boldsymbol{SO(3)}\) and Lifts to the 3-Sphere}\label{sec:so3}
In this section and the next, we generalize~\cite{sal} and compute the expected (Riemannian) distance between two random points in a partially oriented flag manifold. We will discuss this computation for all flags for the case \(n=3\). The strategy for each of these computations is similar. The idea is to lift \(SO(3)\) to the unit 3-sphere and carry out computations upstairs. Lifting to \(\Sph^3\) gives us a natural coordinate system for \(SO(3)\) which is then used to describe the Haar measure on \(SO(3)\). The push-forward measure is then used to find an invariant measure on each of the flags, and amounts to only changing the region of integration in each case. Additionally, we'll discuss how we can do these calculations via Monte Carlo integration.

%-----------------------------------------------------------------------------------------------------------------------------------
\subsection{Lifting for Analytic Computation}\label{sec:so3 distance}
We begin with a known calculation~\cite{sal} of the expected distance between random points in \(SO(3)\). The calculations for partial flags will follow from a refinement of this approach. First we'll briefly describe the Haar measure on \(SO(3)\) as parametrized by \(\Sph^3\). It is not hard to see that \(SO(3)\) is double covered by the unit 3-sphere. Indeed, \(\Sph^3\) has a natural group structure as the unit quaternions, which act on the quaternions by conjugation. This action clearly fixes the real line since the reals commute with quaternions under multiplication. Hence the purely imaginary quaternions are also invariant under the action of unit quaternions. The restriction of the action to the purely imaginary quaternions is an isometry. To see this, suppose that \(q\) is a unit quaternion and that \(x\) and \(y\) are both purely imaginary quaternions. Then
\[
|qxq^{-1} - qyq^{-1}| = |q||x-y||q^{-1}| = |x-y|.
\]

By identifying the purely imaginary quaternions with $\R^3$, this calculation shows that $\Sph^3$ acts by isometries on \(\R^3\). In fact, this action is by rotations: if \(u,n \in \R^3\) are unit vectors, which we interpret as purely imaginary quaternions, and we define \(q = \cos\theta + \sin\theta n \in \Sph^3 \), then a straightforward calculation shows that
\begin{equation}\label{eq:s3 rotation}
quq^{-1} = (\cos\theta +\sin\theta n )u(\cos\theta - \sin\theta n ) = \cos 2\theta\, u + \sin 2 \theta (n \times u) + (1-\cos 2\theta) (u \cdot n) n,
\end{equation}
which is the Rodrigues formula for rotation of \(u\) around the axis \(n\) by an angle \(2\theta\). Therefore, $\Sph^3$ acts on $\R^3$ by rotation; moreover, since \(q\) and \(-q\) both induce the same rotations, it readily follows that \(SO(3)\) is double covered by \(\Sph^3\).

The double covering of \(SO(3)\) provides us a natural way to parametrize \(SO(3)\): each element of $SO(3)$ corresponds to an antipodal pair of points in $\Sph^3$, so each rotation corresponds to a (almost always unique) point in the hemisphere of \(\Sph^3\) where the first coordinate is positive. We can use hyperspherical coordinates on $\Sph^3$ to parametrize \(SO(3)\):
\begin{align*}
x &= \cos \varphi_1\\
y &= \sin \varphi_1 \cos \varphi_2\\
z &= \sin \varphi_1 \sin \varphi_2 \cos \varphi_3\\
w &= \sin \varphi_1 \sin \varphi_2 \sin \varphi_3, 
\end{align*}
where \(\varphi_1  \in [0,\pi/2],~\varphi_2 \in [0,\pi]\) and \(\varphi_3 \in [0,2\pi)\). The volume form on $\Sph^3$ is 
\[
	\sin^2 \varphi_1 \sin\varphi_2 ~d\varphi_1 \wedge d\varphi_2 \wedge d\varphi_3,
\]
but we need to make a slight adjustment to write down the volume form on \(SO(3)\). If \(0 \leq \theta \leq \pi/2\) and \(n\) is a purely imaginary unit quaternion, then the point $q = \cos \theta + \sin \theta n$ has positive first coordinate and lies at a distance \(\theta\) from \(1\) in \(\Sph^3\). However, looking at~\eqref{eq:s3 rotation}, the corresponding element of \(SO(3)\) is at Frobenius distance $2\theta$ from the identity. Therefore, the map \(\Sph^3 \to SO(3)\) scales distances by 2 and hence 3-dimensional volumes by \(2^3=8\), so we conclude that, with respect to hyperspherical coordinates,
\[
	\dVol_{SO(3)} = 8 \sin^2 \varphi_1 \sin\varphi_2 ~d\varphi_1 \wedge d\varphi_2 \wedge d\varphi_3.
\]

 With this in mind, an easy calculation shows that
\[
\vol(SO(3)) = \int_0^{2\pi}\!\!\! \int_0^\pi\!\! \int_0^{\pi/2}\!\! 8 \sin^2 \varphi_1 \sin\varphi_2 d\varphi_1 d\varphi_2 d\varphi_3 = 8\pi^2,
\]
as expected since the first column lies on the unit 2-sphere, with an area of \(4\pi\), the second column lies on a perpendicular unit circle of circumference \(2\pi\), and the third column is determined by the first two.

Now, we can compute the expected distance between two random points in \(SO(3)\). Since the Frobenius distance is equivariant with respect to the left action of \(SO(3)\) on itself, we can rotate one of the two points to the identity and, equivalently, compute the expected distance from a random point in \(SO(3)\) to the identity. Equivalently, since an element of \(SO(3)\) is a rotation by an angle \(\theta\) around an axis \(n\), and \(\theta\) is exactly the Frobenius distance to the identity, we are computing the expected rotation angle of a random element of \(SO(3)\).\footnote{This can also be computed by integrating the rotation angle against the known density of the circular real ensemble~\cite{Girko:1985ju}.}

In order to compute the expected rotation, we note that if \(q = \cos\varphi_1 + n \sin\varphi_1\), the angle of rotation is \(2\varphi_1\). In other words, up in \(\Sph^3\) we are computing the expectation of twice the distance $\varphi_1$ from a random point in the hemisphere of points with positive real part to the identity element \(1\):
\begin{align*}
\E[d;SO(3)] & = \frac{1}{\vol(SO(3))}\!\! \int \limits_{SO(3)} \!\!\! 2\varphi_1 \dVol_{SO(3)}\\
& = \frac{2}{\pi^2}  \int_0^{2\pi} \!\!\! \int_0^\pi \!\! \int_0^{\pi/2} \!\!\! \varphi_1 \sin^2 \varphi_1 \sin\varphi_2 \, d\varphi_1 d\varphi_2 d\varphi_3= \frac{2}{\pi} + \frac{\pi}{2}.
\end{align*}
This calculation agrees with that found in \cite{sal}.

\subsection{Monte Carlo Methods}

This computation can also be done numerically using Monte Carlo integration. First, we will describe how to randomly generate matrices in \(SO(n)\). With this, we will be able to easily compute expected distances on \(SO(3)\). From here, we will be able to use the same algorithms with slight modifications for similar computations on flag and partially oriented flag manifolds. 

Following Chikuse~\cite{Chikuse:2003eg}, we can generate random orthogonal matrices in \Cref{alg:RandSO} by applying Gram--Schmidt to a random Gaussian matrix.

\begin{algorithm}
  \caption{Random Special Orthogonal Matrix}\label{alg:RandSO}
  \begin{algorithmic}[1]
    \Function{RandSO}{$n$}
    \State \(A \gets \text{random $n\times n$ Gaussian}\)
    \State \(Q \gets \textproc{GramSchmidt}{(A)}\)
    \If {\(\det(Q) = 1\)}\\
    	\Return {\(Q\)}
	\Else{ \(Q \gets E_{1,2} Q\)} \Comment{\(E_{1,2}\) is the elementary row swap matrix}\\
	\Return {\(Q\)}
    \EndIf
    \EndFunction
  \end{algorithmic}
\end{algorithm}

The next algorithm depends on the geodesic distance on \(SO(n)\). Fix \(A,B \in SO(n)\), and let \(AB^T = U\Lambda U^*\) be the spectral decomposition of \(AB^T\). Then because \(B\) is an isometry,  
\[
d(A,B) = d(AB^T,I) = d(U\Lambda U^*,I) = d(U \Lambda, U) = d(\Lambda,I).
\] This means the geodesic distance between \(A\) and \(B\) depends only on the eigenvalues of \(AB^T\). Thus when approximating the integral via a Monte Carlo algorithm, we only need to generate one random special orthogonal matrix at each step, which reduces computation time significantly. 

In fact, the distance \(d(A,I)\) is described in \cite{edelman} and is given by
\[
d(A,I) = \left(\frac12 \sum_{k = 1}^n |\log \mu_k|^2\right)^{1/2}.
\]
where \(\mu_k\) are the eigenvalues of \(A\). The factor of \(1/2\) ensures that we aren't overusing \(\arg \mu_k\), since the eigenvalues come in complex conjugate pairs. \Cref{alg:expso} describes the Monte Carlo experiment. Using Algorithms~\ref{alg:RandSO} and~\ref{alg:expso} with \(n=3\) and \(N = 10,\!000,\!000\) gives the approximation \({\E[d;SO(3)] \approx 2.207478465}\), with absolute error
\[
\left| \frac{2}{\pi} + \frac{\pi}{2} - 2.207478465\right| \approx 0.000062365.
\]

\begin{algorithm}
\caption{Calculating Expected Distance}\label{alg:expso}
\begin{algorithmic}[1]
\State \(D \gets [0]*N\) \Comment{Begin with a list of \(N\) zeroes}
\For{$i\gets 1 , N$}
   \State \(A \gets \textproc{RandSO}(n)\)
   \State \(D(i) \gets d(A,I)\)
\EndFor
\end{algorithmic}
\Return \(\textproc{mean}(D)\)
\end{algorithm}

%----------------------------------------------------------------------------------------------------------------------------------
%----------------------------------------------------------------------------------------------------------------------------------
\section{Expected Distances Between Partially Oriented Flags}\label{sec:expected distances}

Following the example computations on \(SO(3)\), we now compute the expected distance between two random points on each (partially oriented) flag manifold \(\Fl(\lambda;P)\) obtained from \(SO(3)\). In all but one case we'll be able to find an analytic expression for the expectation, and in all cases we'll get a numerical result from a Monte Carlo experiment. 

Notice that for \(\lambda\) of length \(k\), \(\sigma \in S_k\), and \(P\) a partition of \(I_\lambda\), we have \(\Fl(\lambda;P) \cong \Fl(\sigma \cdot \lambda; \sigma \cdot P)\). This is true in general, but in the special case of \(SO(3)\) we also observe that if partitions \(P\) and \(P'\) of \(I_\lambda\) have the same order, then \(\Fl(\lambda;P)\) and \(\Fl(\lambda;P')\) are diffeomorphic. Combined with the previous fact, we see that \(\Fl(\lambda;P) \cong \Fl(\sigma \cdot \lambda;\widehat{P})\) so long as the partitions \(P\) of \(I_\lambda\) and \(\widehat{P}\) of \(I_{\sigma \cdot \lambda}\) are the same size. This reduces the number of computations to be done.

In each case, we can use a trick already mentioned above to simplify our calculations. Since each \(\Fl (\lambda; P) = SO(3)/SG_\lambda^P\) is a homogeneous space, and since we will always choose the Riemannian submersion metric on the quotient which is invariant under the left \(SO(3)\) action, we can always move one point to the orbit of the identity, and hence the expected distance between two random points in \(\Fl (\lambda; P)\) is the same as the expected distance between a single random point and the identity orbit (or any other preferred point).

%----------------------------------------------------------------------------------------------------------------------------------
\subsection{\(\Sph^2\) and \(\R P^2\)}

In this section we'll consider the simplest (oriented) flag manifolds derived from \(SO(3)\).\footnote{For completeness, we can consider the trivial flag manifold \(SO(3)/SO(3)\) which clearly has expected distance 0.} Suppose that \(\lambda = (1,2)\). Let $P_C=\{\{1\},\{2\}\}$ and $P_T=\{\{1,2\}\}$ denote the complete and trivial partitions of \(I_\lambda\), respectively. Then we have the identifications \(\Sph^2 \cong \Fl(\lambda; P_C)\) and \(\mathbb{RP}^2 \cong \Fl(\lambda; P_T)\). 

Using standard spherical coordinates, the area form on \(\Sph^2\) is \(\dArea_{\Sph^2} = \sin \varphi \, d\varphi \wedge d\theta\), where \(\varphi \in [0,\pi]\) is the polar angle and \(\theta \in [0,2\pi]\) is the azimuthal angle. Indeed, we verify from the volume formula in \Cref{prop:flag volumes} that $\vol(\Fl((1,2);\{\{1\},\{2\}\}))=4\pi$, the area of the unit sphere. Since \(S^2\) is homogeneous, it suffices to consider the distance between a point and the north pole, which is simply the polar angle \(\varphi\). Hence we have 
\[
\E[d;\Fl((1,2);P_C)] =  \frac{1}{\Area(\Sph^2)} \int \limits_{\Sph^2} \varphi \dArea_{\Sph^2} = \frac{1}{4\pi} \int_0^{2\pi}\!\!\! \int_0^\pi \!\! \varphi \sin \varphi\, d\varphi d\theta = \frac{\pi}{2}.
\]

This is what we expect, since a point is just as likely to be in the northern hemisphere as in the southern. The calculation can also be carried out using Monte Carlo integration: generate a 3-vector \(v\) with independent, normally distributed entries. Normalize \(v\) to get \(\hat{v}\) which is uniformly distributed on the sphere. As above, it suffices to compute the angle between \(\hat{v}\) and the north pole. Taking the average of \(N\) samples gives the numerical estimate; indeed, an experiment with \(N=10,\!000,\!000\) yields the estimate \(\E[d;\Fl((1,2);P_C)] \approx 1.570989\), with an absolute error of \(\left|\frac{\pi}{2} - 1.570989\right|\approx 0.000193\).

The story for \(\mathbb{RP}^2\) is similar. The map \(\pi: \Sph^2 \to \mathbb{RP}^2\) given by \(p \mapsto [p] := \{p,-p\}\) is a Riemannian submersion, so the area form on \(\mathbb{RP}^2\) is \(\pi_* \dArea_{\Sph^2} = \sin \varphi \,d\varphi \wedge d\theta\), and computing the distance between points in $\mathbb{RP}^2$ is equivalent to computing the minimum distance between two pairs of antipodal points in $\Sph^2$, so our computation reduces to determining the expected angle between a random point in the northern hemisphere and the north pole:
\[
\E[d;\Fl((1,2);P_T)] = \frac{1}{\Area(\mathbb{RP}^2)} \int \limits_{\mathbb{RP}^2} \varphi \dArea_{\mathbb{RP}^2} = \frac{1}{2\pi} \int_0^{2\pi}\!\!\!\! \int_0^{\pi/2} \varphi\sin \varphi\, d\varphi d\theta = 1.
\]

To carry this out with Monte Carlo methods, take \(\hat{v}\) as before, but compute the average of \({\arccos|\hat{v} \cdot e_3|}\) where \(e_3\) is the point at the north pole. An experiment with \(N=10,\!000,\!000\) produces the estimate \(\E[d;\Fl((1,2);P_T)] \approx 0.999975\).

%------------------------------------------------------------------------------------------------------------------------------------
\subsection{Flags with \(\lambda = (1,1,1)\)}
In this section we consider flags with the ordered partition \(\lambda = (1,1,1)\). First we'll describe the problem of finding expected distances via Monte Carlo integration. Following this, we'll take a look at each individual case. Suppose that \(P\) is any partition of \(I_\lambda\). In this case \(SG_\lambda^P\) contains only a finite number of elements in \(SO(3)\) as we have already seen in \Cref{fin}. Computing the expected distance on \(SO(3)/SG_\lambda^P\) is hardly any different from what we've done before: we simply take the minimum distance from each point in the orbit \(gSG_\lambda^P\) to the identity as in \Cref{alg:exphom}.

\begin{algorithm}
\caption{Expected Distance on $SO(3)/SG_\lambda^P$}\label{alg:exphom}
\begin{algorithmic}[1]
\State \(D \gets [0]*N\) \Comment{Begin with a list of \(N\) zeroes}
\For{$i\gets 1 , N$}
   \State \(A \gets \textproc{RandSO}(3)\)
   \State \(\text{orbitdata} \gets [0]*|SG_\lambda^P|\) \Comment{\(|SG_\lambda^P|\) denotes the order of \(SG_\lambda^P\)}
   \For{$ j \gets 1, |SG_\lambda^P|$}
   \State \( \text{orbitdata}(j) \gets d(Ah_j,I)\) \Comment{Each \(h_j\) denotes a distinct element of \(SG_\lambda^P\)}
   \EndFor
   \State \(D(i) \gets \textproc{min}(\text{orbitdata})\)
\EndFor
\end{algorithmic}
\Return \(\textproc{mean}(D)\)
\end{algorithm}

%----------------------------------------------------------------------------------------------------------------------------------
%----------------------------------------------------------------------------------------------------------------------------------
\subsubsection{Oriented Flags}

When \(P = P_C\) is the complete partition, the oriented flag manifold \(\Fl(\lambda; P) \cong SO(3)\) and we saw in \Cref{sec:so3 distance} that the expected distance between two random points in this space is 
\[
	\E[d;\Fl((1,1,1);P_C)] = E[d;SO(3)] = \frac{2}{\pi} + \frac{\pi}{2}.
\]

%----------------------------------------------------------------------------------------------------------------------------------
%----------------------------------------------------------------------------------------------------------------------------------
\subsubsection{Partially Oriented Flags}
For the case of the partially oriented flags we may choose \(P\) to be \(P_1,P_2,\) or \(P_3\) as in \Cref{fin}. In any case, \(SG_\lambda^P\) will contain two matrices. Without loss of generality, suppose that \(P = P_1\). We can readily apply \Cref{alg:exphom} with  \(N = 10,\!000,\!000 \) to see that the expected distance is \(\approx 1.78548266\).

To get an analytic expression, the strategy is to lift the computation to \(\Sph^3\) as in \Cref{sec:so3}. Since the composite map \(\Sph^3 \to SO(3) \to SO(3)/SG_\lambda^P\) is (up to the scale factor 2) a Riemannian submersion, the distance between a pair of points in \(\Fl(\lambda; P) = SO(3)/SG_\lambda^P\) is equal to twice the minimum distance between the sets of preimages in \(\Sph^3\). 

As usual, the expected distance between two random points in \(\Fl(\lambda; P)\) is the same as the expected distance between a single random point and any preferred point, which we will take to be the identity coset \(I SG_\lambda^P\). In turn, the identity coset consists of the identity and the rotation by angle \(\pi\) around the axis \((1,0,0)\), so its preimages in \(\Sph^3\) are \(\{1,-1,i,-i\}\); in general, if \(g \in SO(3)\) is rotation by \(\theta\) around an axis \(n\), then the preimage of \(g SG_\lambda^P\) in \(\Sph^3\) is \(\{q,-q,iq,-iq\}\), where \(q = \cos \frac{\theta}{2} + \sin \frac{\theta}{2} n\). In turn, since multiplication by \(-1\) and multiplication by \(i\) are isometries, 
\[
	d(\{1,-1,i,-i\},\{q,-q,iq,-iq\}) = d(1,\{q,-q,iq,-iq\}),
\]
where the minimum distance will be achieved by the element of \(\{q,-q,iq,-iq\}\) which is closer to \(1\) than to \(-1\), \(i\), or \(-i\).

Therefore, the expected distance between a random coset \(g SG_\lambda^P\) and the identity coset (and hence, between two random points in \(\Fl(\lambda; P)\)), is simply the expectation of twice the distance from \(1\) to a random point in the subset of \(\Sph^3\) which is closer to \(1\) than to \(-1\), \(i\), or \(-i\). In terms of Cartesian coordinates \((x,y,z,w)\), this is precisely the set \(\{(x,y,z,w): x \geq |y|\}\); in hyperspherical coordinates, \(x \geq |y|\) is equivalent to \(0 \leq \phi_1 \leq \arctan(\sec \phi_2)\). 

Hence, the volume of the partial flag is given by
\begin{multline*}
\vol(\Fl(\lambda; P))= \!\!\! \int \limits_{\Fl(\lambda; P)} \!\!\! \dVol_{\Fl(\lambda; P)}=\int_0^{2\pi}\!\!\! \int_0^{\pi}\!\! \int_0^{\arctan(\sec \varphi_2)}\!\! 8\sin^2 \varphi_1 \sin\varphi_2 ~d\varphi_1d\varphi_2d\varphi_3 \\
= 2\int_0^{2\pi}\!\!\! \int_0^{\pi/2}\!\!\!\! \int_0^{\arctan(\sec \varphi_2)}\!\! 8\sin^2 \varphi_1 \sin\varphi_2 ~d\varphi_1d\varphi_2d\varphi_3,
\end{multline*}
where we use symmetry to restrict to \(\varphi_2 \in [0,\pi/2)\) and avoid any issues with \(\arctan(\sec \varphi_2)\).

\emph{Mathematica} will happily compute the above integral to be \(4\pi^2\), which agrees with the volume calculation from \Cref{prop:flag volumes}, but the limits of integration on the innermost integral make this unpleasant to compute by hand.

\begin{figure}[htbp]
	\centering
		\includegraphics[height=3in]{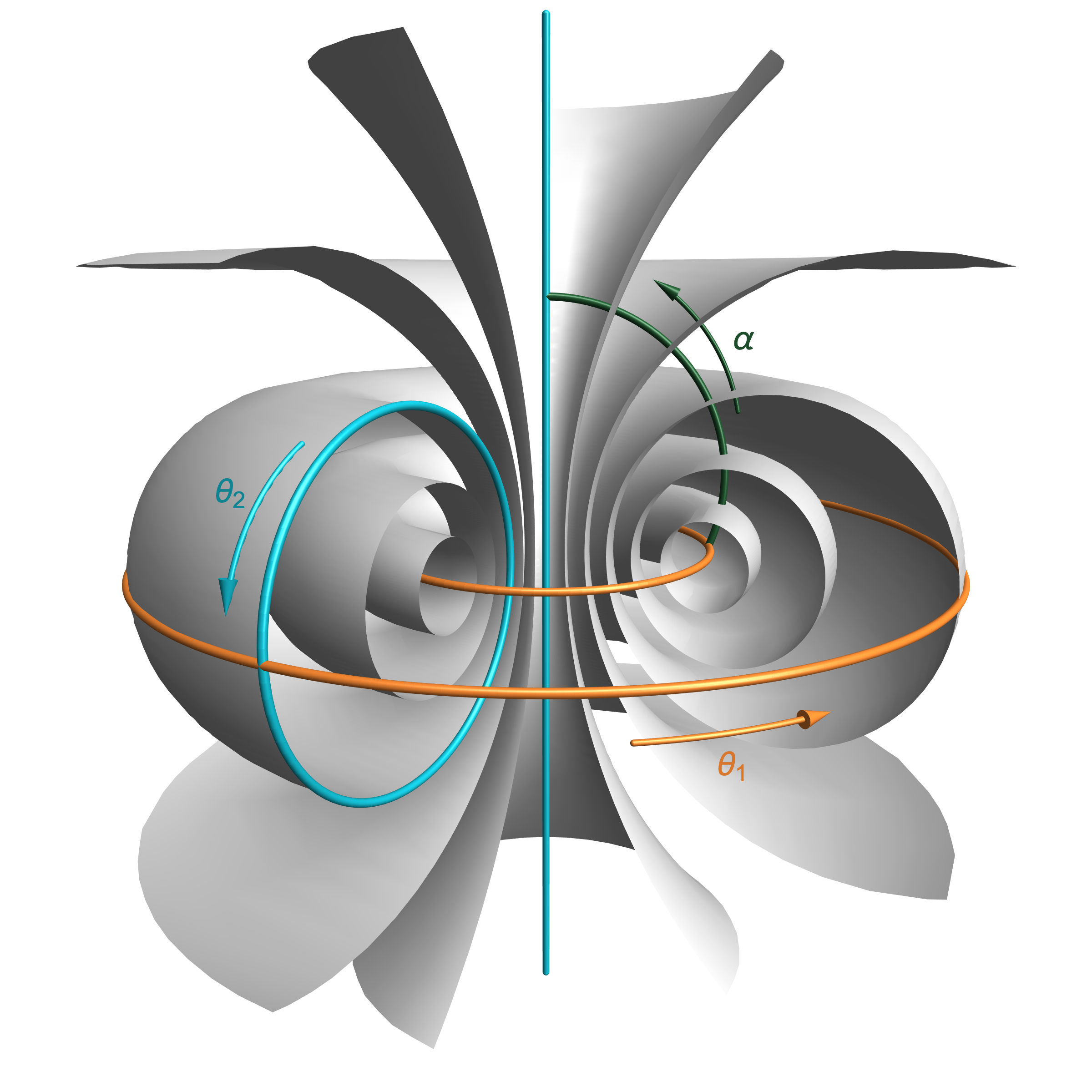}
	\caption{Stereographic projection to \(\R^3\) of join coordinates on \(\Sph^3\). \(\theta_1\) and \(\theta_2\) are simply the arguments of the two complex coordinates, while \(\alpha\) gives the angle a given vector makes with the \(z_1\)-plane. The level sets of \(\alpha\) are generically tori, collapsing to the unit circle in the \(z_1\)-plane when \(\alpha = 0\) and to the unit circle in the \(z_2\)-plane (which stereographically projects to the \(z\)-axis) when \(\alpha = \pi/2\).}
	\label{fig:join}
\end{figure}

We can make things easier by changing coordinate systems. Thinking of \(\Sph^3\) as the unit sphere in \(\C^2\), we will interpret points on \(\Sph^3\) as pairs \((z_1,z_2)\in \C^2\) with \(|z_1|^2 + |z_2|^2 = 1\). Then the condition $x > |y|$ is equivalent to \(|\arg z_1| < \frac{\pi}{4}\) -- which will be a lens-shaped region in \(\Sph^3\); the analogous region on \(\Sph^2\) is a lune -- so it is desirable to use $\arg z_1$ as one of our coordinates. In fact, \(|z_1|^2 + |z_2|^2 = 1\) means that
\[
	(z_1, z_2) = (\cos \alpha\, e^{i \theta_1}\!, \sin \alpha\, e^{i \theta_2})
\]
for \(\alpha \in [0,\pi/2]\) and \(\theta_1,\theta_2 \in [0,2\pi]\). Notice that each value of \(\alpha\) determines a torus in \(\Sph^3\) except the extremal values \(\alpha=0\) and \(\alpha=\pi/2\), where the torus collapses to the unit circle in the \(z_1\)- or \(z_2\)-plane. See \Cref{fig:join} for a visualization.

We can write these coordinates -- which we call \textbf{join coordinates} because they give a concrete realization of \(\Sph^3\) as the topological join of two copies of \(\Sph^1\) -- in terms of Cartesian coordinates as
\begin{align*}
x &= \cos \alpha \cos \theta_1 \\
y &= \cos \alpha \sin \theta_1 \\
z &= \sin \alpha \cos \theta_2 \\
w &= \sin \alpha \sin \theta_2
\end{align*}
and the volume form on \(\Sph^3\) is easily computed to be \(\dVol_{\Sph^3} = \cos \alpha \sin \alpha \,d\alpha \wedge d\theta_1 \wedge d\theta_2\). In these coordinates, the volume of \(\Fl(\lambda; P)\) is easy to compute by hand, as it reduces to a simple \(u\)-substitution:
\[
	\vol(\Fl(\lambda; P)) = \int_0^{2\pi} \!\!\!\int_{-\pi/4}^{\pi/4} \int_0^{\pi/2}\!\!\!  8	\cos \alpha \sin \alpha \, d\alpha\, d\theta_1 d\theta_2 = 4\pi^2
\]

We want to compute the expectation of twice the spherical distance from \(1\) to a random point in the region \(-\frac{\pi}{4} < \theta_1 < \frac{\pi}{4}\). The distance is simply \(\arccos(\cos \alpha \cos \theta_1)\), and so the expected distance between two random points in \(\Fl(\lambda; P)\) is
\begin{multline*}
\frac{1}{\vol(\Fl(\lambda; P))} \int_{-\pi}^\pi \int_{-\pi/4}^{\pi/4} \int_0^{\pi/2} \!\! 2\arccos(\cos \alpha \cos \theta_1)  8 \cos \alpha \sin \alpha \,d\alpha\, d\theta_1 d\theta_2\\
= \frac{4}{\pi^2} \int_{-\pi}^\pi \int_{-\pi/4}^{\pi/4} \int_0^{\pi/2} \!\! \arccos(\cos \alpha \cos \theta_1)  \cos \alpha \sin \alpha \,d\alpha\, d\theta_1 d\theta_2 \\
= \frac{16}{\pi} \int_0^{\pi/4} \!\!\! \int_0^{\pi/2} \!\! \arccos(\cos \alpha \cos \theta_1)  \cos \alpha \sin \alpha \,d\alpha\, d\theta_1
\end{multline*}
by integrating out \(\theta_2\) and using the fact that the integrand is even in \(\theta_1\). 

This is slightly tedious but essentially straightforward to compute using the substitutions \(u = \cos \alpha\) and \({\sin v = u \cos \theta_1}\), producing the first half of \Cref{thm:main theorem}, which we restate in slightly more general form:
\begin{theorem}\label{thm:partial flag}
	For $i=1,2,3$, the expected distance between two random points in \(\Fl((1,1,1);P_i)\) is
	\[
		\E[d;\Fl((1,1,1);P_i)] = 1 + \frac{\pi}{4}.
	\]
\end{theorem}

Note that 
\[
\left| 1 + \frac\pi4 - 1.78548266\right| \approx 0.0000844966
\]
which once again shows that our Monte Carlo experiment does a good job of approximating the analytic result.
%----------------------------------------------------------------------------------------------------------------------------------
%----------------------------------------------------------------------------------------------------------------------------------
\subsubsection{Flags}
Finally, we consider the case when \(P\) is the trivial partition, so \(SG_\lambda^P\) contains four matrices as in \Cref{fin} and \(\Fl(\lambda; P) \cong \Fl(1,2,3)\), the manifold of (complete) flags on \(\R^3\). We may readily apply \Cref{alg:exphom} with \(N = 10,\!000,\!000\) to find that the expected distance between two random points in \(\Fl(1,2,3)\) is
\(\approx 1.311751687\). 

Finding an explicit integral for this calculation only requires a few additional changes to the bounds from the previous section. The orbit of the identity lifts to \(\{\pm 1, \pm i, \pm j, \pm k\}\), which are the vertices of the regular 16-cell dual to the standard hypercube. Other orbits are simply rotated copies of the orbit of the identity, so they also get lifted to regular 16-cells. Hence, computing the expected distance between points in \(\Fl(1,2,3)\) is equivalent to computing the expected distance from $1$ to all the other unit quaternions which are closer to $1$ than any of the other integer points on \(\Sph^3\). 

This set is nothing but the radial projection of the cube 
\[
	\{(x,y,z,w) : x=1, -1\leq |y|,|z|,|w| \leq 1\}
\] 
to \(\Sph^3\), namely those points with \(x \geq |y|\), \(x \geq |z|\), and \(x \geq |w|\). This spherical cube is partitioned into \(3! \times 2^3\) spherical tetrahedra, each congruent to the spherical tetrahedron $x \geq y \geq z \geq w \geq 0$. To integrate over this tetrahedron in hyperspherical coordinates we get the limits of integration \(\varphi_1 \in [0, \arctan(\sec \varphi_2)]\), \(\varphi_2 \in [0,\arctan(\sec \varphi_3)]\), and \(\varphi_3 \in [0,\pi/4]\). Since the simplex is \(\frac{1}{48}\) of the full spherical cube, we see that
\[
\vol(\Fl(1,2,3)) = 48 \int_0^{\pi/4}\!\!\!\! \int_0^{\arctan(\sec(\varphi_3))}\!\!\!\! \int_0^{\arctan(\sec(\varphi_2))}\!\! 8\sin^2 \varphi_1 \sin \varphi_2 \,d\varphi_1 d\varphi_2d\varphi_3;
\]
using \emph{Mathematica}'s numerical integration algorithm reassuringly yields a number numerically indistinguishable from \(2\pi^2\), which we know from \Cref{prop:flag volumes} is the true volume.

Therefore, the expected distance between two points on \(\Fl(1,2,3)\) is given by the integral
\[
\frac{48}{2\pi^2}\int_0^{\pi/4}\!\!\!\! \int_0^{\arctan(\sec(\varphi_3))}\!\!\!\! \int_0^{\arctan(\sec(\varphi_2))}\!\! 2\varphi_1 8 \sin^2 \varphi_1 \sin \varphi_2 \,d\varphi_1 d\varphi_2d\varphi_3.
\]
In join coordinates, the corresponding integral turns out to be
\[
	\frac{4}{2\pi^2} \int_{-\pi/4}^{\pi/4} \int_{-\pi/4}^{\pi/4} \int_0^{\arctan \left( \frac{\cos \theta_1}{\cos \theta_2}\right)}\!\! 2 \arccos(\cos \alpha \cos \theta_1)\, 8 \cos \alpha \sin \alpha\, d\alpha\, d\theta_1 d\theta_2.
\]

It is not clear how to evaluate either of these integrals exactly, though the former can be reduced to the one-dimensional integral given in the second half of \Cref{thm:main theorem}, which we restate:

\begin{theorem}\label{thm:full flag}
	The expected distance between two random \((1,2,3)\)-flags is
	\begin{multline*}
		\E[d;\Fl(1,2,3)] = \frac{3\pi}{2} + \frac{96}{\pi^2} \bigintsss_0^{\pi/4}\! \left[ \arctan\left( \tan^2 \left(\frac{\arctan(\sec \varphi_3)}{2}\right)\right) - \frac{\arctan^2 \left(\sqrt{1+\sec^2 \varphi_3}\right)}{\sqrt{1+\sec^2\varphi_3}}\right]\! d\varphi_3.
	\end{multline*}
\end{theorem}

Of course, this can be numerically evaluated to an arbitrary degree of precision; to 20 digits, we can compute that the expected distance between two random points in \(\Fl(1,2,3)\) is $1.3117250347224445929$. Using the PSLQ algorithm~\cite{Ferguson:1999gn}, this does not appear to be in the vector space over \(\Q\) generated by \(1, \pi,\) and \(\frac{1}{\pi}\), unlike all the other expected values we have computed.

%----------------------------------------------------------------------------------------------------------------------------------
%---------------------------------------------------------------------------------------------------------------------------------
\section{Conclusion and Open Questions}

The manifolds of partially oriented flags described in this paper have the potential to be a natural home for data which is interpretable in terms of nested subspaces, of which some may be oriented. We encourage mathematicians, engineers, and other practitioners who are beginning to see the usefulness of flag manifolds in data analysis to keep these spaces in mind.

While partially oriented flag manifolds have been studied by algebraic topologists~\cite{Lam:1975kj, Sankaran:1987wh,Sankaran:1997cj}, we have not been able to find much evidence of their being studied by geometers. Given the fundamental role played by flag varieties in algebraic geometry, it would be interesting to find a more algebraic definition or interpretation of these spaces. We also look forward to further geometric computation on these spaces, whether it be computing expected distances on a larger class of partially oriented flag manifolds or determining other geometric quantities of interest.

For data living in a metric space, the expected distance between two random points in that space is a simple statistical baseline for distance-based similarity measures. While this quantity can often be estimated effectively using Monte Carlo techniques, it would be interesting to see analytic expressions for expected distance on spaces which are nice geometrically or useful as data models.

\section*{Acknowledgments}

This project grew out of conversations in the Pattern Analysis Lab at Colorado State University, and we are very grateful to all of the participants in the lab for ongoing inspiration. This work was partially supported by the National Science Foundation (ATD \#1712788, CCF--BSF:CIF \#1830676, CP) and the Simons Foundation (\#354225, CS).

\bibliographystyle{plain}
\bibliography{bib}
\end{document}